\documentclass{article}
\usepackage{arxiv}
\usepackage{hyperref}
\usepackage[utf8]{inputenc} % allow utf-8 input
\usepackage[T1]{fontenc}    % use 8-bit T1 fonts
\usepackage{url}            % simple URL typesetting
\usepackage{booktabs}       % professional-quality tables
\usepackage{amsfonts}       % blackboard math symbols
\usepackage{nicefrac}       % compact symbols for 1/2, etc.
\usepackage{microtype}      % microtypography
\usepackage{lipsum}		% Can be removed after putting your text content
\usepackage{graphicx}
\usepackage{doi}

\usepackage{indentfirst}
\usepackage{float}
\usepackage{tikz}
\usepackage{enumitem}           % enumerate com i) ii)
\usepackage{amsmath,amsbsy,amsfonts,amssymb}   
\usepackage{times}
\usepackage{mathptmx} 
\usepackage{subfigure}

\title{Optimal vaccination strategies in the control of an infectious disease: a SEIRV model for administration of two vaccines}

\author{ \href{https://orcid.org/0009-0007-8132-7118}{\includegraphics[scale=0.06]{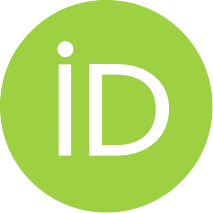}\hspace{1mm}Nelson L. Santos Junior} \\
	Department of Mathematics\\
	Federal University of Pernambuco\\
	Recife, Pernambuco, Brazil \\
	\texttt{nelson.leal@ufpe.br} \\
	\And
	\href{https://orcid.org/0000-0003-4855-4136}{\includegraphics[scale=0.06]{orcid.pdf}\hspace{1mm}João A. M. Gondim} \\
	Department of Mathematics\\
	Federal University of Pernambuco\\
	Recife, Pernambuco, Brazil \\
	\texttt{joao.gondim@ufpe.br} \\
}

\renewcommand{\shorttitle}{\textit{arXiv} Template}

\hypersetup{
pdftitle={A template for the arxiv style},
pdfsubject={q-bio.NC, q-bio.QM},
pdfauthor={David S.~Hippocampus, Elias D.~Striatum},
pdfkeywords={First keyword, Second keyword, More},
}

\begin{document}
\maketitle

%-------------------
\begin{abstract}
In this paper, we study the optimal control for an SEIR model adapted to the vaccination strategy of susceptible individuals. There are factors associated with a vaccination campaign that make this strategy not only a public health issue but also an economic one. In this case, optimal control is important as it minimizes implementation costs. We consider the availability of two vaccines with different efficacy levels, and the control indicates when each vaccine should be used. The optimal strategy specifies in all cases how vaccine purchases should be distributed. For similar efficacy values, we perform a sensitivity analysis on parameters that depend on the intrinsic characteristics of the vaccines. Additionally, we investigate the behavior of the number of infections under the optimal vaccination strategy.
\end{abstract}
%-------------------

%-------------------
%-------------------

% ----------------------------------------------------------
\section{Introduction}
% ----------------------------------------------------------

In recent decades, mathematical models have become an important tool for investigating and predicting the dynamics of disease spread. These studies help us understand how different containment strategies should work \cite{mathematicalmodels}. Many epidemiological models have been pivotal in recent years, and the demand for developing such models was particularly high during the COVID-19 pandemic, an event that promises to be a historic milestone of the 21st century \cite{Horton}.

Vaccination is a containment strategy that significantly reduces the risk of infection, severe cases, hospitalization, and death, providing not only individual protection to vaccinated people, but also community-level protection \cite{importancevaccine}. It is estimated that in the case of COVID-19, a global reduction of $63\%$ in total deaths was achieved during the first year of vaccination. In low-income countries, it is further estimated that the number of avoided deaths would increase by $111\%$ if the vaccine coverage target of $40\%$, as established by the World Health Organization (WHO), had been met by each country by the end of 2021 \cite{impactovacinacao1}.

Historically, many studies have focused on mathematical models of vaccination in the context of various diseases, such as hepatitis and measles \cite{vacinacaohepatite}, influenza \cite{vacinacaoinfluenza}, tuberculosis \cite{vacinacaotuberculose}, and malaria \cite{malaria}. Recently, numerous models have been developed to investigate different aspects of vaccination against SARS-CoV-2 during the COVID-19 pandemic, including vaccine distribution \cite{alocacao1, alocacao2}, impact on the number of deaths \cite{impactovacinacao1, impactovacinacao2}, sensitivity analysis of model parameters \cite{sensibilidade}, the behavior of the basic reproductive number \cite{mainreference}, and other related contexts.

Disease containment strategies generate financial costs for the government. For vaccination campaigns, for example, some of these costs are: purchase, storage, transportation, and distribution of vaccines, and installation of vaccination sites and medical staffing \cite{impact3}, in addition to investments in public awareness campaigns. Thus, a vaccination campaign is not only a public health measure but also a strategy with significant economic implications, and numerous studies have sought to investigate this impact \cite{impact3, impact4, impact5}. In this context, it is essential to develop mathematical models that examine containment strategies while managing associated economic costs, thereby providing insights to support government decisions in implementing these policies.

In epidemiological models, an optimal control problem can be used to define the ideal value of an important parameter of the model, so that it maximizes or minimizes a certain other variable present; then we can use optimal control to minimize a financial cost associated with a disease containment campaign. The theory of optimal control has been applied to disease models such as: HIV \cite{hiv1,hiv2,hiv3}, tuberculosis \cite{tub1,tub2}, influenza \cite{gripe}, and general epidemics \cite{behncke,mateus}. Recently, many studies have used optimal control for COVID-19 study \cite{demasse,shen,obsu,zamir}. Gondim and Machado \cite{Gondim2} used optimal control to analyze a SEIR model with quarantine of susceptible people which represents a total lockdown system during the COVID-19 pandemic. Santos Junior and Gondim \cite{junior} used optimal control to analyze a model for quarantine only of infected people which is closer to the proper definition of quarantine.

The goal of this paper is work with the optimal control on a SEIR model adapted for vaccination campaign of susceptible people with two types of vaccines, we will calculate numerically these controls for describe the ideal use a given type of vaccine, minimizing not only the number of infected people, but also the economic costs associated with the use of this vaccine. We hope that the results can suggest government decisions when implementing the vaccination campaigns. The results indicates when the government must use a vaccine or other depending of your efficacy, in this cases, we indicates the purchase percentage of each vaccine. When the two vaccines have efficacy values very near, we present a study of sensibility for those parameters that depend of other characteristic of vaccines. Finally, we analyze the number of infected in three contexts: when we only use the first vaccine, when we use both vaccines according optimal control and when we only use the second vaccine.

% ----------------------------------------------------------
\section{The epidemiology model}
% ----------------------------------------------------------

We will consider a modified SEIR model with a one-dose vaccination strategy for susceptible individuals. We know that in the case of COVID-19, for example, a vaccine with partial protection was sufficient to control the SARS-CoV-2 pandemic \cite{Sarah, Monia}, so we will work with vaccines that are not $100\%$ effective. Here, we consider two vaccines, $V_1$ and $V_2$, with efficacy $\theta_1$ and $\theta_2$, respectively, where $\theta_2 < \theta_1 < 1$.

From \cite{Tentori}, the efficacy $\theta_i$, for $i=1,2$, is given by $\theta_i = (P - P_i)/P$, where $P$ and $P_i$ indicate the attack rates for unvaccinated and vaccinated individuals, respectively, i.e., the probability that a contact results in transmission. Thus, we have $P_i = P(1-\theta_i)$ for $i=1,2$. From \cite{Martcheva}, the transmission rate in the SEIR model is $\beta = cP$, where $c$ is the contact rate. Therefore, for vaccinated individuals, the transmission rate is given by $\beta_i = cP_i = cP(1-\theta_i) = \beta (1-\theta_i)$, meaning the vaccine efficacy directly reduces the transmission rate of the disease.

The efficacy values chosen for use in the simulations are motivated by real COVID-19 vaccine efficacies. For example, the Pfizer vaccine is $91\%$ effective \cite{pfizer}, the AstraZeneca vaccine is $74\%$ effective \cite{astrazeneca}, the Janssen vaccine is $67\%$ effective \cite{janssen}, and the Coronavac vaccine is $51\%$ effective \cite{coronavac}. These references are based on phase 3 trials in various countries with different population groups; some efficacies are guaranteed for two doses of the vaccine, but for simplicity, we assume these efficacies are valid for one dose of the vaccine in our study.

Let $S(t), V_1(t), V_2(t), E(t), I(t)$, and $R(t)$ be the numbers of susceptible individuals, individuals vaccinated with vaccine $V_1$, individuals vaccinated with vaccine $V_2$, exposed individuals, infected individuals, and recovered individuals, respectively. We assume that the total population 
\begin{equation*}
N(t)=S(t)+V_1(t)+V_2(t)+E(t)+I(t)+R(t)
\end{equation*}
is constant, as we are considering only a short time frame in comparison to the demographic time scale.

Vaccination was associated with high short-term protection against SARS-CoV-2, but this protection decreases considerably after a few months \cite{hall}. This strategy helps contain disease spread but does not completely eliminate the risk of infection, as reinfections post-vaccination have been reported \cite{reinfeccao1, reinfeccao2}. Therefore, our model includes the effect of a potential reinfection by the disease.

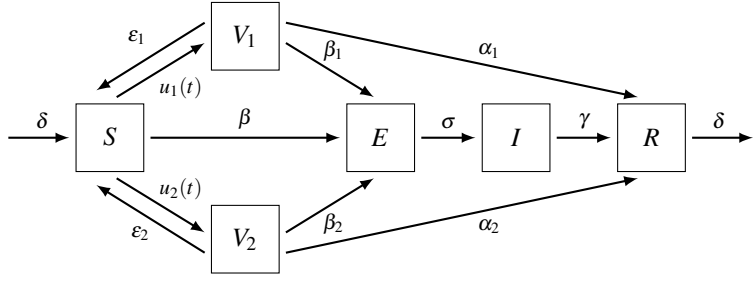
\begin{figure}[H]
\label{diagram}
\centering
\begin{tikzpicture}[scale=0.45]
 %%%%%%%%%%%%%%% S %%%%%%%%%%%%%%% 

\draw[black] (-9, -1) rectangle (-7,1) {}; 
\node[circle,inner sep=0.7pt,label=above:{$S$}] at (-8,-0.6) {};

\draw[thick][-latex] (-5.2,3.4) -- (-8.4,1.4);
\node[circle,inner sep=0.7pt,label=above:{ \scalebox{0.8}{$\epsilon_1$}}] at (-7.2, 2.4) {};

\draw[thick][-latex] (-5.2,-3.4) -- (-8.4,-1.4);
\node[circle,inner sep=0.7pt,label=above:{ \scalebox{0.8}{$\epsilon_2$}}] at (-7.1, -3.4) {};

\draw[thick][-latex] (-11,0) -- (-9.2,0);
\node[circle,inner sep=0.7pt,label=above:{\scalebox{0.8}{$\delta$}}] at (-10,-0.1) {};

 %%%%%%%%%%%%%%% V1 %%%%%%%%%%%%%%% 

\draw[black] (-5, 2) rectangle (-3,4) {}; 
\node[circle,inner sep=0.7pt,label=above:{$V_1$}] at (-4,2.3) {};

\draw[thick][-latex] (-7.8,1.2) -- (-5.2,2.8);
\node[circle,inner sep=0.7pt,label=above:{ \scalebox{0.8}{$u_1(t)$}}] at (-5.9, 0.8) {};

 %%%%%%%%%%%%%%% V2 %%%%%%%%%%%%%%% 

\draw[black] (-5, -2) rectangle (-3,-4) {}; 
\node[circle,inner sep=0.7pt,label=above:{$V_2$}] at (-4,-3.7) {};

\draw[thick][-latex] (-7.8,-1.2) -- (-5.2,-2.8);
\node[circle,inner sep=0.7pt,label=above:{ \scalebox{0.8}{$u_2(t)$}}] at (-5.9, -2.2) {};

 %%%%%%%%%%%%%%% E %%%%%%%%%%%%%%% 

\draw[black] (-1,-1) rectangle (1,1); 
\node[circle,inner sep=0.7pt,label=above:{$E$}] at (0,-0.6) {};

\draw[thick][-latex] (-6.8,0) -- (-1.2,0);
\node[circle,inner sep=0.7pt,label=above:{ \scalebox{0.8}{$\beta$}}] at (-4, -0.1) {};

\draw[thick][-latex] (-2.8,2.8) -- (-0.2,1.2);
\node[circle,inner sep=0.7pt,label=above:{ \scalebox{0.8}{$\beta_{1}$}}] at (-1.4, 1.9) {};

\draw[thick][-latex] (-2.8,-2.8) -- (-0.2,-1.2);
\node[circle,inner sep=0.7pt,label=above:{ \scalebox{0.8}{$\beta_{2}$}}] at (-1.4,-3.2) {};

 %%%%%%%%%%%%%%% I %%%%%%%%%%%%%%% 
 
\draw[black] (3, 0-1) rectangle (5,1);
\node[circle,inner sep=0.7pt,label=above:{$I$}] at (4,-0.6) {};

\draw[thick][-latex] (1.2,0) -- (2.8,0);
\node[circle,inner sep=0.7pt,label=above:{ \scalebox{0.8}{$\sigma$}}] at (2,-0.1) {};

 %%%%%%%%%%%%%%% R %%%%%%%%%%%%%%% 
 
\draw[black] (7,-1) rectangle (9,1) {}; 
\node[circle,inner sep=0.7pt,label=above:{$R$}] at (8,-0.6) {};

\draw[thick][-latex] (5.2,0) -- (6.8,0);
\node[circle,inner sep=0.7pt,label=above:{\scalebox{0.8}{$\gamma$}}] at (6,-0.1) {};

\draw[thick][-latex] (-2.8,3.4) -- (7.6,1.2);
\node[circle,inner sep=0.7pt,label=above:{ \scalebox{0.8}{$\alpha_{1}$}}] at (3.2, 2) {};

\draw[thick][-latex] (-2.8,-3.4) -- (7.6,-1.2);
\node[circle,inner sep=0.7pt,label=above:{ \scalebox{0.8}{$\alpha_{2}$}}] at (3.2,-3.2) {};

\draw[thick][-latex] (9.2,0) -- (11,0);
\node[circle,inner sep=0.7pt,label=above:{\scalebox{0.8}{$\delta$}}] at (10,-0.1) {};

\end{tikzpicture}
\vspace*{0.3cm}\caption{Compartment diagram of the model.}
\end{figure}

Figure \ref{diagram} illustrates the dynamics between the compartments of our model.

The susceptible population decreases when individuals move into the exposed class at a transmission rate $\beta$, are vaccinated with vaccine $V_1$ at a control rate $u_1(t)$, or are vaccinated with vaccine $V_2$ at a control rate $u_2(t)$. The controls represent the fraction of individuals vaccinated per unit time at $t$. Here, our boundary condition for the control is 
\begin{equation}\label{limitcontrol}
0 \leq u(t) \leq 1,
\end{equation}
since there is no maximum time to be vaccinated, and the minimum time is 1 day.

In the model, there are three possibilities for a vaccinated individual:
\begin{enumerate}
\item Leave the class of those vaccinated with immunity and move into the class of recovered people at a rate $\alpha_i$, for $i=1,2$;
\item Being vulnerable after leaving the vaccinated class with possibility of contracting the disease, moving into the susceptible class to a rate $\epsilon_i$, for $i=1,2$;
\item Contract the disease before leaving the vaccinated class becoming a exposed individual and transmitting the disease at a rate $\beta_i=1-\theta_i$, for $i=1,2$;
\end{enumerate}

The exposed population decreases when individuals move into the infected class at a latency rate $\sigma$. The recovered population increases at an exit rate $\gamma$ and decreases at a reinfection rate $\delta$. Note that in this model, people can recover by acquiring vaccine-induced immunity without necessarily experiencing an infection.

The system of equations is
\begin{equation}\label{model}
\left\{\begin{array}{l}
\displaystyle \frac{d S}{dt} = - \frac{\beta S I}{N} - u_1(t) S - u_2(t) S + \epsilon_1 V_1 + \epsilon_2 V_2 + \delta R \\
\displaystyle \frac{d V_1}{dt} = u_1(t) S - \frac{\beta_1 V_1 I}{N} - \epsilon_1 V_1 - \alpha_1 V_1\\
\displaystyle \frac{d V_2}{dt} = u_2(t) S - \frac{\beta_2 V_2 I}{N} -\epsilon_2 V_2 - \alpha_2 V_2 \\
\displaystyle \frac{d E}{dt} = \frac{\beta S I}{N} + \frac{\beta _1 V_1 I}{N} + \frac{\beta_2 V_2 I}{N} - \sigma E  \\
\displaystyle \frac{d I}{dt} = \sigma E - \gamma I \\
\displaystyle \frac{d R}{dt} = \gamma I - \delta R + \alpha_1 V_1 + \alpha_2 V_2
\end{array}\right.
\end{equation}

All parameters are nonnegative and were taken from \cite{mainreference}, where the Bayesian Markov Chain Monte Carlo (MCMC) method \cite{MCMC} shows the convergence of the parameters. In that reference, only one vaccine is considered, providing a single value for the immunity rate and for the rate of return to the susceptible class. In our study, we incorporate two vaccines, and we need two values for each parameter. For simplicity, we initially assume $\alpha_1=\alpha_2$ and $\epsilon_1=\epsilon_2$. The parameter values for our numerical simulations are described in Table \ref{parameters}.

\begin{table}[H]
\caption{\label{parameters} Parameters values.}
\begin{center}
\begin{tabular}{c|c|c} 
\hline
\textbf {Parameter} & \textbf{Description} & \textbf {Value} \\
\hline
$\beta$ & transmission rate & 0.45  \\
\hline
$\sigma$ & latency rate & 0.25  \\
\hline
$\gamma$ & recovered rate & 0.07 \\ 
\hline
$\delta$ & reinfection rate & 0.65 \\
\hline
$\alpha_1$ & $V_1$ vaccine immunity rate & 0.08 \\
\hline
$\alpha_2$ & $V_2$ vaccine immunity rate & 0.08  \\
\hline
$\epsilon_1$ & rate of return to susceptible class of $V_1$ vaccine & 0.54 \\
\hline
$\epsilon_2$ & rate of return to susceptible class of $V_2$ vaccine & 0.54 \\
\hline
\end{tabular}
\end{center}
\end{table}

The number of susceptible individuals is considered equal to the total population, as the proportion of exposed, infected, and recovered individuals is very small in comparison to the total population, as in \cite{Gondim2}. We assume that there are no vaccinated individuals at the start of the simulation. Our initial time point uses data from Brazil as of May 8, 2020 \cite{Gondim2, Worldmeters}. The number of recovered individuals is calculated as the difference between the total number of cases and the number of deaths. The mean incubation period of the disease is estimated to be around 5 days \cite{Lauer}, so the number of exposed individuals is determined by subtracting the number of active cases on May 8 from the number of active cases on May 13. The initial conditions are described in Table \ref{initialconditions}.

\begin{table}[H]
\caption{\label{initialconditions} Description of initial conditions.}
\begin{center}
\begin{tabular}{c|c} 
\hline
\textbf {Group} & \textbf{Initial Condition} \\
\hline
Susceptible & 200000000 \\
\hline
Vaccinated with $V_1$ vaccine & 0 \\		     
\hline
Vaccinated with $V_2$ vaccine & 0 \\		     
\hline
Exposed & 65124 \\
\hline
Infected & 76603 \\		     
\hline
Recovery & 65124 \\		     
\hline
\end{tabular}
\end{center}
\end{table}

% ----------------------------------------------------------
\section{The optimization problem}
% ----------------------------------------------------------

Following system \eqref{model}, we minimize the functional 
\begin{equation}\label{functional}
J = \int_0^T  \big( I(t) + B_1 u_1^2(t) + B_2 u_2^2(t) \big)~dt. 
\end{equation}

In the formula \eqref{functional}, $T$ is the duration of vaccination and the parameters $B_1$ and $B_2$ are the costs of implementation of $V_1$ and $V_2$ vaccines, respectively. We assume that $B_1, B_2>0$, as well as that 
\begin{equation}
B = B_1 + B_2,   
\end{equation}
where $B$ is the total control cost.

Sufficient conditions for the existence of optimal controls follow from standard results from optimal control theory, we can use Theorem 2.1 in Joshi et al. \cite{Joshi} to show that the optimal control exists.

Pontryagin's maximum principle \cite{Pontryagin,Lenhart} establish that optimal controls are solutions of the Hamiltonian system with Hamiltonian function
\begin{equation}\label{hamiltonian}
\begin{aligned}
\displaystyle H &= I + B_1 u_1^2 + B_2 u_2^2 + u_1 S \big(\lambda_{V_1} - \lambda_{S}\big) + u_2 S \big(\lambda_{V_2} - \lambda_{S}\big) \\
&+ \frac{\beta S I}{N} \big(\lambda_E - \lambda_S\big) + \frac{\beta_1 V_1 I}{N} \big(\lambda_E - \lambda_{V_1}\big) + \frac{\beta_2 V_2 I}{N} \big(\lambda_E - \lambda_{V_2}\big) \\
&+ \epsilon_1 V_1 \big(\lambda_S - \lambda_{V_1}\big) + \epsilon_2 V_2 \big(\lambda_S - \lambda_{V_2}\big) + \alpha_1 V_1 \big(\lambda_R - \lambda_{V_1}\big) \\
&+ \alpha_2 V_2 \big(\lambda_R - \lambda_{V_2}\big) + \sigma E \big(\lambda_I - \lambda_E\big) + \delta R \big(\lambda_{S} - \lambda_R\big) + \gamma I \big(\lambda_{R} - \lambda_I\big),
\end{aligned}
\end{equation} 
involving the state variables $S, V_1, V_2, E, I, R$ and the adjoint variables $\lambda_S$, $\lambda_{V_1}$, $\lambda_{V_2}$,$\lambda_E$ $\lambda_I$, $\lambda_R$, which satisfy the adjoint system \eqref{adjointsystem}. 
\begin{equation}\label{adjointsystem}
\left\{\begin{array}{cl}
\displaystyle \frac{d \lambda_S}{dt} & = \displaystyle  \frac{I}{N^2} \bigg[ \beta_1 V_1 \big( \lambda_E -\lambda_{V_1} \big) + \beta_2 V_2 \big( \lambda_E -\lambda_{V_2} \big) + \beta \big(N - S \big) \big( \lambda_S -\lambda_E \big) \bigg] \\
 & + ~ u_1 \big( \lambda_{S} - \lambda_{V_1} \big) + u_2 \big( \lambda_{S} -\lambda_{V_2} \big) \\
 
\displaystyle \frac{d \lambda_{V_1}}{dt} & = \displaystyle  \frac{I}{N^2} \bigg[ \beta S \big( \lambda_E -\lambda_S \big) + \beta_2 V_2 \big( \lambda_E -\lambda_{V_2} \big) + \beta_1 \big(N - V_1\big) \big( \lambda_{V_1} -\lambda_E \big) \bigg] \\
& + ~\epsilon_1 \big( \lambda_{V_1} - \lambda_{S} \big) + \alpha_1 \big( \lambda_{V_1} -\lambda_{R} \big) \\

\displaystyle \frac{d \lambda_{V_2}}{dt} & = \displaystyle  \frac{I}{N^2} \bigg[ \beta S \big( \lambda_E -\lambda_S \big) + \beta_1 V_1 \big( \lambda_E -\lambda_{V_1} \big) + \beta_2 \big(N - V_2\big) \big( \lambda_{V_2} -\lambda_E \big) \bigg] \\
& + ~\epsilon_2 \big( \lambda_{V_2} - \lambda_{S} \big) + \alpha_2 \big( \lambda_{V_2} -\lambda_{R} \big) \\

\displaystyle \frac{d \lambda_{E}}{dt} & = \displaystyle  \frac{I}{N^2} \bigg[ \beta S \big( \lambda_E -\lambda_S \big) + \beta_1 V_1 \big( \lambda_E -\lambda_{V_1} \big) + \beta_2 V_2 \big( \lambda_{E} -\lambda_{V_2} \big) \bigg] \\
&+~\sigma \big( \lambda_{E} - \lambda_{I} \big)  \\

\displaystyle \frac{d \lambda_I}{dt} & = \displaystyle  \frac{N-I}{N^2} \bigg[ \beta S \big( \lambda_S -\lambda_{E} \big) + \beta_1 V_1 \big( \lambda_{V_1} -\lambda_{E} \big) +   \beta_2 V_2 \big( \lambda_{V_2} -\lambda_E \big) \bigg]  \\
&+~\gamma \big( \lambda_{I} - \lambda_{R} \big) - 1 \\

\displaystyle \frac{d \lambda_{R}}{dt} & = \displaystyle  \frac{I}{N^2} \bigg[ \beta S \big( \lambda_E -\lambda_S \big) + \beta_1 V_1 \big( \lambda_E -\lambda_{V_1} \big) + \beta_2 V_2 \big( \lambda_{E} -\lambda_{V_2} \big) \bigg] \\
&+~\delta \big( \lambda_{R} - \lambda_{S} \big) 
\end{array}\right.
\end{equation}

The adjoint variables also must satisfy the transversality conditions
\begin{equation}\label{transversalitycondition}
\lambda_S(T) = \lambda_{V_1}(T) = \lambda_{V_2}(T) = \lambda_E(T) = \lambda_I(T) = \lambda_R(T) = 0.
\end{equation}

Finally, the optimality conditions come from solving
\begin{equation*}
\frac{\partial H}{\partial u}\Bigg|_{u=u^*}=0, 
\end{equation*}
this results in
\begin{equation}
\displaystyle u_i^*=\frac{S\big(\lambda_{S}-\lambda_{V_i}\big)}{2 B_i}, \quad i=1,2. 
\end{equation}

Since we are considering bounded controls, by \eqref{limitcontrol} the optimal control $u^*$ are calculated \cite{Lenhart} using
\begin{equation}\label{optimalycondition}
\displaystyle u^*=\min \left\{1, \max \left\{0, \frac{S\big(\lambda_{S}-\lambda_{V_i}\big)}{2 B_i}\right\}\right\}, \quad i=1,2.
\end{equation}

Uniqueness of the optimal controls (at least for small enough T) also follow from standard results, such as Theorem 2.3 in Joshi et al. \cite{Joshi} . 

The numerical solutions of systems \eqref{model} and \eqref{adjointsystem} can be found using the forward-backward sweep method \cite{Lenhart}. We have a initial condition in $t=0$ for system \eqref{model} and a initial condition in $t=T$ for system \eqref{adjointsystem}, the algorithm of sweep uses these conditions and works with the following steps:
    
\textbf{Step 1}: Starts with an initial guess of the controls $u_1$ and $u_2$, solves \eqref{model} using the initial condition for state variable and the values for $u_i$ forward in time;

\textbf{Step 2}: Uses initial condition \eqref{transversalitycondition} and the results provided in Step 1 for the state variables and for the optimal control, to solve \eqref{adjointsystem} backward in time, the new controls are defined by \eqref{optimalycondition};

\textbf{Step 3}: This process continues until a convergence is obtained.

% ----------------------------------------------------------
\section{Discussion}
% ----------------------------------------------------------

Our optimal control problem depends on the implementation costs $B_1$ and $B_2$ of the two vaccines $V_1$ and $V_2$, respectively. Here, the costs for the vaccines will be distributed as $B_i = \theta_i \cdot 10^4$ for $i=1,2$; i.e., the implementation cost of a vaccine is proportional to its efficacy power. This decision was made due to the difficulty of finding specific references on government spending for vaccine procurement.

Next, see the plots the graphs of the optimal controls $u_1(t)$ and $u_2(t)$ for vaccination durations of 60, 120, and 180 days. In all cases, $u_1$ represents the implementation of the vaccine with higher efficacy, therefore, higher cost.

\begin{figure}[H]
\begin{center}
\includegraphics[scale=0.9, trim=20mm 115mm 20mm 115mm]{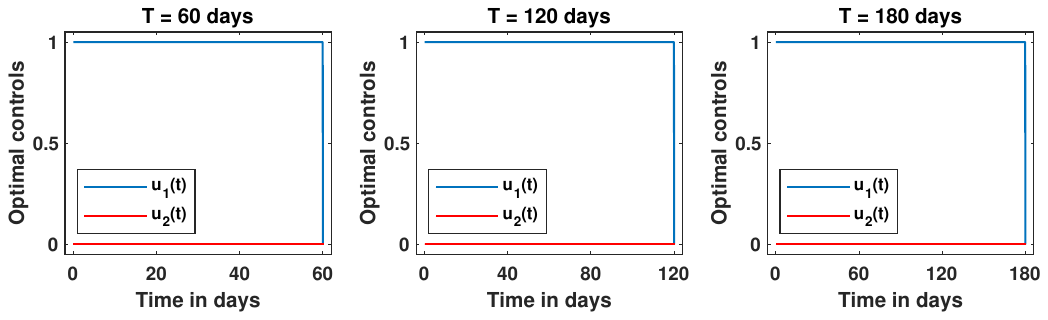}
\end{center}
\caption{\label{control1} The optimal controls for use of vaccines with effectiveness of $91\%$ and $51\%$ for different vaccination campaigns lengths.}
\end{figure}

\begin{figure}[H]
\begin{center}
\includegraphics[scale=0.9, trim=20mm 115mm 20mm 115mm]{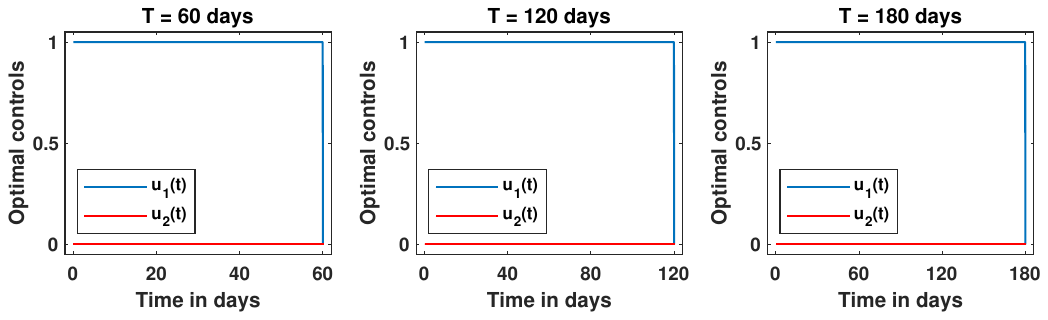}
\end{center}
\caption{\label{control2} The optimal controls for use of vaccines with effectiveness of $74\%$ and $51\%$ for different vaccination campaigns lengths.}
\end{figure}

\begin{figure}[H]
\begin{center}
\includegraphics[scale=0.9, trim=20mm 115mm 20mm 115mm]{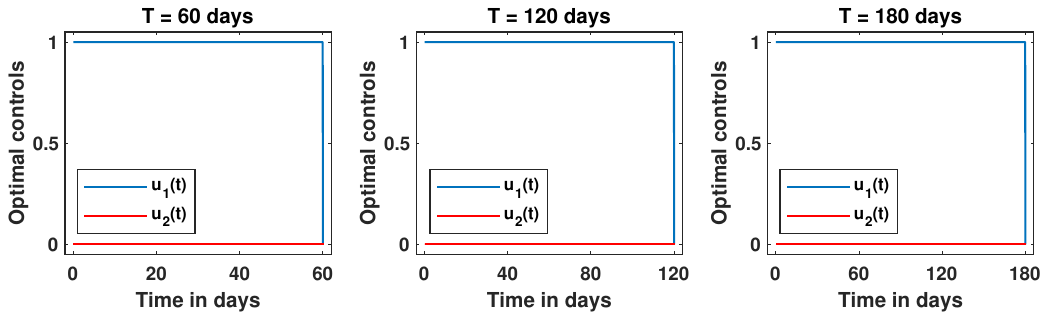}
\end{center}
\caption{\label{control3} The optimal controls for use of vaccines with effectiveness of $67\%$ and $51\%$ for different vaccination campaigns lengths.}
\end{figure}

\begin{figure}[H]
\begin{center}
\includegraphics[scale=0.9, trim=20mm 115mm 20mm 115mm]{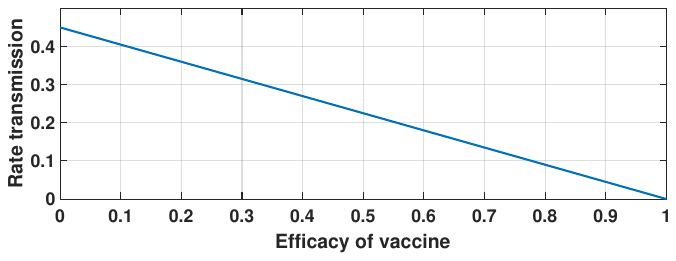}
\end{center}
\caption{\label{ratetransmission} Behavior of rate transmission for vaccinated people.}
\end{figure}

In Figures \ref{control1}, \ref{control2}, and \ref{control3}, regardless of the higher economic cost, the optimal control suggests prioritizing the use of the $V_1$ vaccine with an efficacy of $91\%$, $74\%$, or $67\%$ over the $V_2$ vaccine, which has an efficacy of $51\%$. This preference is likely due to the impact of the vaccine’s efficacy on reducing the disease transmission rate. Specifically, for vaccinated individuals, the transmission rate decreases linearly, as shown in Figure \ref{ratetransmission}. In this context, note that $\beta_2 = 0.2378$, while smaller than the usual rate for unvaccinated individuals, $\beta=0.45$, remains significantly higher than $\beta_1 = 0.1601$, $\beta_1 = 0.1262$, or $\beta_1= 0.0436$.

\begin{figure}[H]
\begin{center}
\includegraphics[scale=0.9, trim=20mm 115mm 20mm 115mm]{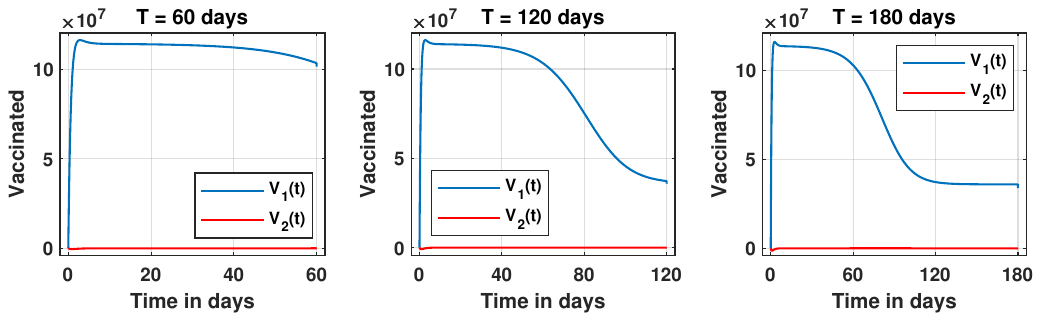}
\end{center}
\caption{\label{vaccinated1} Vaccination curves for the use of vaccines with effectiveness of $91\%$ and $51\%$ for different vaccination campaigns lengths.}
\end{figure}

\begin{figure}[H]
\begin{center}
\includegraphics[scale=0.9, trim=20mm 115mm 20mm 115mm]{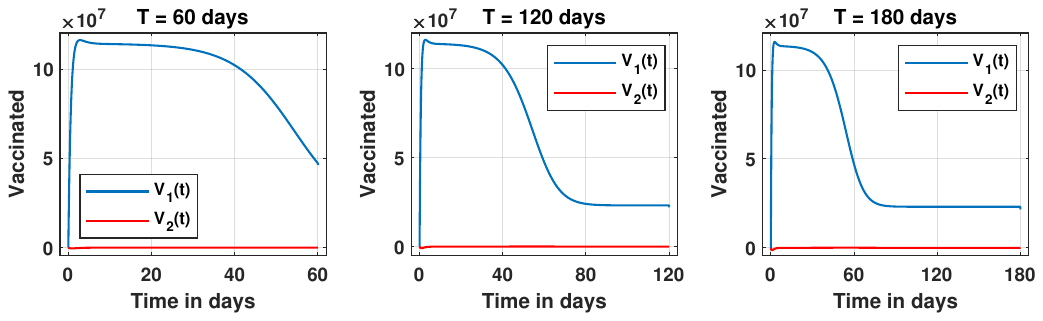}
\end{center}
\caption{\label{vaccinated2} Vaccination curves for the use of vaccines with effectiveness of $74\%$ and $51\%$ for different vaccination campaigns lengths.}
\end{figure}

\begin{figure}[H]
\begin{center}
\includegraphics[scale=0.9, trim=20mm 115mm 20mm 115mm]{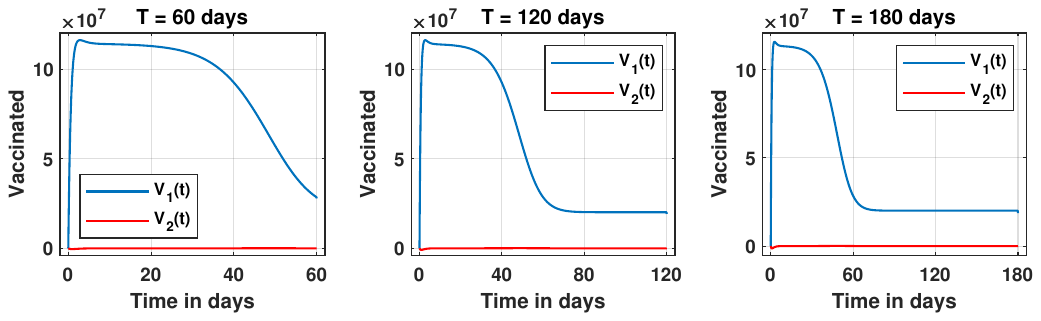}
\end{center}
\caption{\label{vaccinated3} Vaccination curves for the use of vaccines with effectiveness of $67\%$ and $51\%$ for different vaccination campaigns lengths.}
\end{figure}

In Figures \ref{vaccinated1}, \ref{vaccinated2} and \ref{vaccinated3} see the curve of vaccinated for the controls of the Figures \ref{control1}, \ref{control2} and \ref{control3}, respectively. The number of vaccinated people with the $V_2$ vaccine is practically zero when compared to the size of the susceptible population, this is a expected result considering what we saw in the plots of the optimal controls. After this discussion, it is recommended that the government purchase only the most effective vaccine, a vaccine with a efficacy of $51\%$ only be recommended if it is the only one available. 

\begin{figure}[H]
\begin{center}
\includegraphics[scale=0.9, trim=20mm 115mm 20mm 115mm]{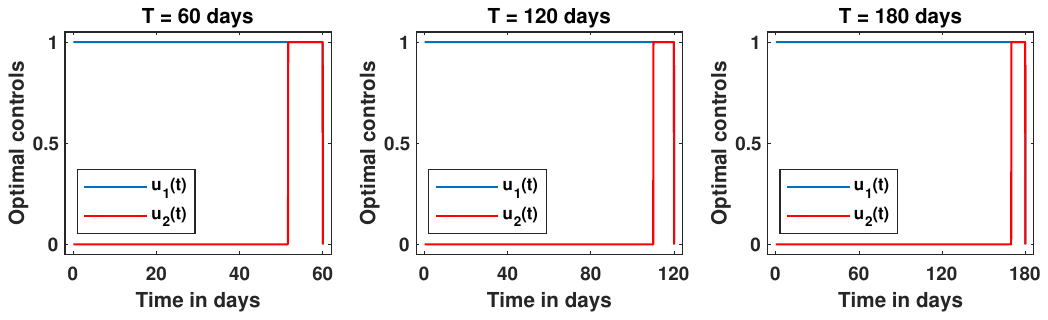}
\end{center}
\caption{\label{control4} The optimal controls for the use of vaccines with effectiveness of $91\%$ and $74\%$ for different vaccination campaigns lengths.}
\end{figure}

\begin{figure}[H]
\begin{center}
\includegraphics[scale=0.9, trim=20mm 115mm 20mm 115mm]{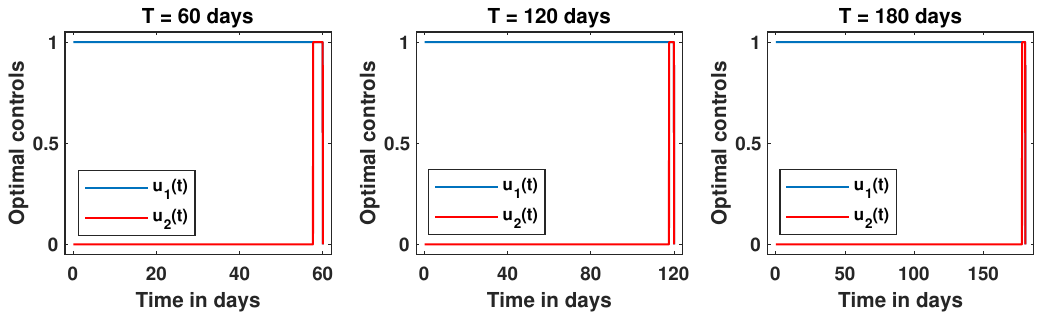}
\end{center}
\caption{\label{control5} The optimal controls for use of vaccines with effectiveness of $91\%$ and $67\%$ for different vaccination campaigns lengths.}
\end{figure}

In Figure \ref{control4}, for most of the campaign duration, the optimal control suggests using the $V_1$ vaccine with $91\%$ efficacy compared to the $V_2$ vaccine with $74\%$ efficacy. However, there is a simultaneous use of both vaccines starting on days 51, 109, and 169 for vaccination campaigns lasting 60, 120, and 180 days, respectively. When the $V_2$ vaccine has $67\%$ efficacy, as shown in Figure \ref{control5}, the optimal control behavior remains almost identical to that in Figure \ref{control4}, with simultaneous use of both vaccines beginning on days 57, 116, and 176 for vaccination campaigns lasting 60, 120, and 180 days, respectively.

We understand that the strategy of the optimal controls in Figures \ref{control4} and \ref{control5} seeks to balance economic expenses, given that the $V_2$ vaccine has a lower implementation cost. Additionally, the lower effectiveness of the $V_2$ vaccine is unlikely to pose issues, as the population has likely already achieved herd immunity through the initial administration of the more effective $V_1$ vaccine.

\begin{figure}[H]
\begin{center}
\includegraphics[scale=0.9, trim=20mm 115mm 20mm 115mm]{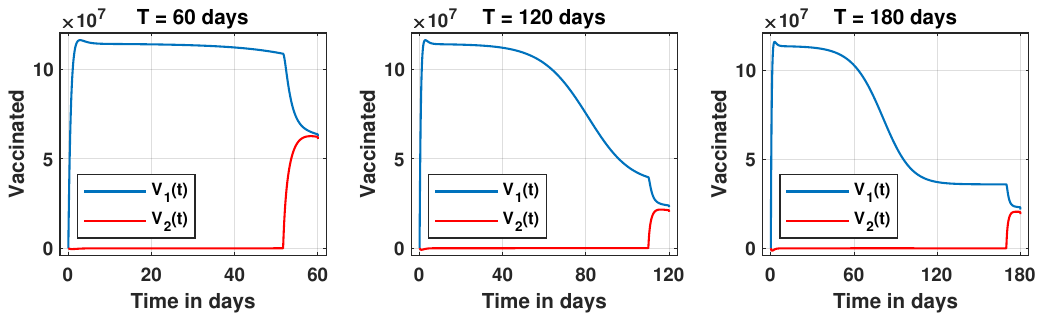}
\end{center}
\caption{\label{vaccinated4} Vaccination curves for the use of vaccines with effectiveness of $91\%$ for $74\%$ for different vaccination campaigns lengths.}
\end{figure}

\begin{figure}[H]
\begin{center}
\includegraphics[scale=0.9, trim=20mm 115mm 20mm 115mm]{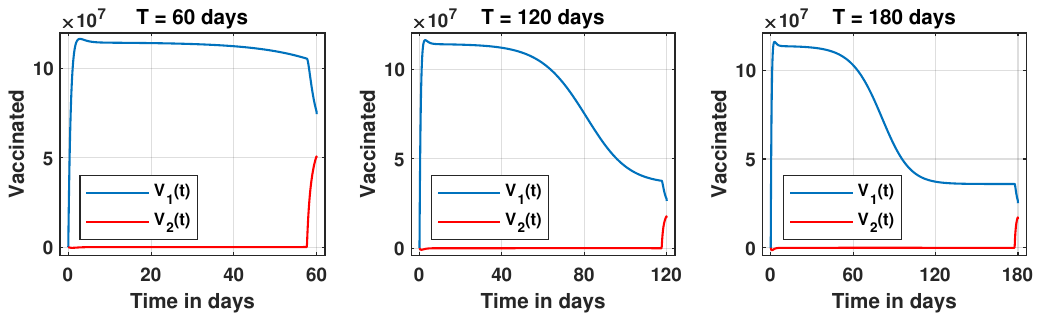}
\end{center}
\caption{\label{vaccinated5} Vaccination curves for the use of vaccines with effectiveness of $91\%$ and $67\%$ for different vaccination campaigns lengths.}
\end{figure}

In Figures \ref{vaccinated4} and \ref{vaccinated5}, see the curve of vaccinated individuals based on the controls shown in Figures \ref{control4} and \ref{control5}. See that due to the use of the $V_2$ vaccine indicated by control, there is an increase in the number of people vaccinated with this vaccine in the final days of the campaign.

Using the Trapezoidal rule for numerical integration \cite{chapra}, we calculated the integral of the vaccination curves $V_1(t)$ and $V_2(t)$ to determine the quantity of vaccines used during the campaign, providing a recommended distribution for government agencies to guide vaccine procurement.

\begin{table}[H]
\caption{\label{purchasedistribution4} Distribution for vaccine procurement according to Figure \ref{vaccinated4}.}
\begin{center}
\begin{tabular}{c|c|c|c}
\hline
\textbf{Vaccine} & \textbf{T = 60 days} & \textbf{T = 120 days} & \textbf{T = 180 days} \\
\hline
$V_1$ & 93,44 \% & 98,18 \% & 98,55 \% \\	     
\hline
$V_2$ &  6,66 \% &  1,82 \% &  1,45 \% \\	     
\hline
\end{tabular}
\end{center}
\end{table}

\begin{table}[H]
\caption{\label{purchasedistribution5} Distribution for vaccine procurement according to Figure \ref{vaccinated5}.}
\begin{center}
\begin{tabular}{c|c|c|c}
\hline
\textbf{Vaccine} & \textbf{T = 60 days} & \textbf{T = 120 days} & \textbf{T = 180 days} \\
\hline
$V_1$ & 98,77 \% & 99,70 \% & 99,76 \% \\	     
\hline
$V_2$ &  1,23 \% &  0,30 \% &  0,24 \% \\	     
\hline
\end{tabular}
\end{center}
\end{table}

\begin{figure}[H]
\begin{center}
\includegraphics[scale=0.9, trim=20mm 115mm 20mm 115mm]{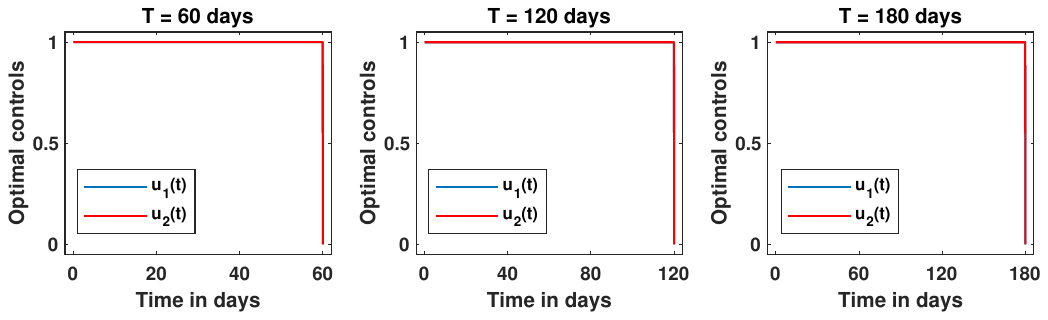}
\end{center}
\caption{\label{control6} The optimal controls for use of vaccines with effectiveness of $74\%$ and $67\%$ and different vaccination campaigns lengths}
\end{figure}

For vaccines with $74\%$ and $67\%$ efficacy, there is a small difference between the efficacy rates $\theta_1,\theta_2$ and the transmission rates $\beta_1,\beta_2$. In this case, as shown in Figure \ref{control6}, the plot of optimal control is overlaid, indicating a simultaneous use of the vaccines throughout the duration of the vaccination campaign.

\begin{figure}[H]
\begin{center}
\includegraphics[scale=0.9, trim=20mm 115mm 20mm 115mm]{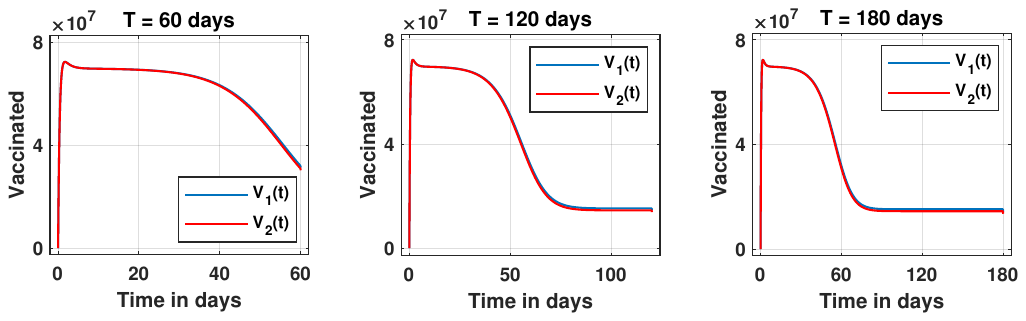}
\end{center}
\caption{\label{vaccinated6} Vaccination curves for the use of vaccines with effectiveness of $74\%$ and $67\%$ and different vaccination campaigns lengths.}
\end{figure}

See in Figure \ref{vaccinated6}, the curve of vaccinated people based on the optimal control in Figure \ref{control6}.

In all cases shown in Figure \ref{vaccinated6}, the number of people vaccinated with the $V_1$ and $V_2$ vaccines remains practically the same throughout the vaccination period. In this context, we also calculate the integral of the vaccinated curves $V_1(t)$ and $V_2(t)$ to observe the quantity of vaccines used during the campaign, as presented in Table \ref{purchasedistribution6}.

\begin{table}[H]
\caption{\label{purchasedistribution6} Distribution for vaccine procurement according to Figure \ref{vaccinated6}.}
\begin{center}
\begin{tabular}{c|c|c|c}
\hline
\textbf{Vaccine} & \textbf{T = 60 days} & \textbf{T = 120 days} & \textbf{T = 180 days} \\
\hline
$V_1$ & 50,11 \% & 50,34 \% & 50,48 \% \\	     
\hline
$V_2$ & 49,89 \% & 49,66 \% & 49,52 \% \\	     
\hline
\end{tabular}
\end{center}
\end{table}

We know that from Figure \ref{control6}, the optimal control does not indicate preference for the use of one vaccine over another with the parameters adopted in Table \ref{parameters}. In this context, we investigated the sensitivity of the parameters related to the immunity rate, the rate of return to the susceptible class, and efficacy, which can vary in the real world according to the characteristics of the vaccines used.

\begin{table}[H]
\caption{\label{sensibility} Values of $\alpha_1$, $\alpha_2$, $\epsilon_1$ and $\epsilon_2$ for sensibility analysis.}
\begin{center}
\resizebox{\textwidth}{!}{%
\begin{tabular}{c|c|c|c}
\hline
$\alpha_1 > \alpha_2 \quad \epsilon_1 > \epsilon_2$ & 
$\alpha_1 < \alpha_2 \quad \epsilon_1 < \epsilon_2$ & 
$\alpha_1 > \alpha_2 \quad \epsilon_1 < \epsilon_2$ &
$\alpha_1 < \alpha_2 \quad \epsilon_1 > \epsilon_2$ \\
\hline
$\alpha_1 = 0.088 \quad \epsilon_1 = 0.594$ & 
$\alpha_1 = 0.08 \quad \epsilon_1 = 0.54$ &
$\alpha_1 = 0.088 \quad \epsilon_1 = 0.54$ &
$\alpha_1 = 0.08 \quad \epsilon_1 = 0.594$ \\
$\alpha_2 = 0.08 \quad \epsilon_2 = 0.54$ & 
$\alpha_2 = 0.008 \quad \epsilon_2 = 0.594$ &
$\alpha_2 = 0.08 \quad \epsilon_2 = 0.594$ &
$\alpha_2 = 0.088 \quad \epsilon_2 = 0.54$ \\
\hline
$\alpha_1 = 0.096 \quad \epsilon_1 = 0.648$ & 
$\alpha_1 = 0.08 \quad \epsilon_1 = 0.54$ &
$\alpha_1 = 0.096 \quad \epsilon_1 = 0.54$ &
$\alpha_1 = 0.08 \quad \epsilon_1 = 0.648$ \\
$\alpha_2 = 0.08 \quad \epsilon_2 = 0.54$ & 
$\alpha_2 = 0.096 \quad \epsilon_2 = 0.648$ &
$\alpha_2 = 0.08 \quad \epsilon_2 = 0.648$ &
$\alpha_2 = 0.096 \quad \epsilon_2 = 0.54$ \\	     
\hline
$\alpha_1 = 0.08 \quad \epsilon_1 = 0.54$ & 
$\alpha_1 = 0.072 \quad \epsilon_1 = 0.486$ &
$\alpha_1 = 0.08 \quad \epsilon_1 = 0.486$ &
$\alpha_1 = 0.072 \quad \epsilon_1 = 0.54$ \\
$\alpha_2 = 0.072 \quad \epsilon_2 = 0.486$ & 
$\alpha_2 = 0.08 \quad \epsilon_2 = 0.54$ &
$\alpha_2 = 0.072 \quad \epsilon_2 = 0.54$ &
$\alpha_2 = 0.08 \quad \epsilon_2 = 0.486$ \\
\hline
$\alpha_1 = 0.08 \quad \epsilon_1 = 0.54$ & 
$\alpha_1 = 0.064 \quad \epsilon_1 = 0.432$ &
$\alpha_1 = 0.08 \quad \epsilon_1 = 0.432$ &
$\alpha_1 = 0.064 \quad \epsilon_1 = 0.54$ \\
$\alpha_2 = 0.064 \quad \epsilon_2 = 0.432$ & 
$\alpha_2 = 0.08 \quad \epsilon_2 = 0.54$ &
$\alpha_2 = 0.064 \quad \epsilon_2 = 0.54$ &
$\alpha_2 = 0.08 \quad \epsilon_2 = 0.432$ \\		     		     
\hline
\end{tabular} }
\end{center}
\end{table}

\begin{figure}[H]
\begin{center}
\includegraphics[scale=0.9, trim=20mm 115mm 20mm 115mm]{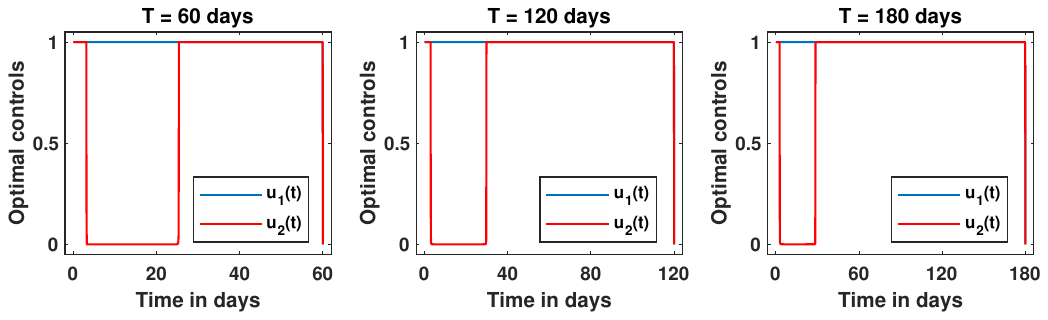}
\end{center}
\caption{\label{sensibility3} Optimal controls for $\alpha_1 = 0.08, \alpha_2 = 0.064, \epsilon_1 = 0.432, \epsilon_2 = 0.54$.}
\end{figure}

In Table \ref{sensibility}, we vary the values of the parameters $\alpha_1$, $\alpha_2$, $\epsilon_1$, and $\epsilon_2$ up and down by $10\%$ and $20\%$, and the optimal controls were calculated for this set of values. For all percentage up variations, the optimal controls exhibit the same behavior as shown in Figure \ref{control6}. For the down variations, there is one case of a change in control, which is illustrated in Figure \ref{sensibility3}.

In Figure \ref{sensibility3}, the reduction in the use of the $V_2$ vaccine is acceptable. Note that we have a reduction of $20\%$ in the parameter values of $\alpha_2$ and $\epsilon_1$. This change in parameters indicates that the $V_2$ vaccine has a smaller recovery power, while the $V_1$ vaccine shows a reduction in the number of people returning to the susceptible class. Furthermore, the $20\%$ variation represents a threshold; that is, a reduction of $19\%$ does not result in a change in the behavior of the optimal controls depicted in Figure \ref{control6}.

\begin{figure}[H]
\begin{center}
\includegraphics[scale=0.9, trim=20mm 115mm 20mm 115mm]{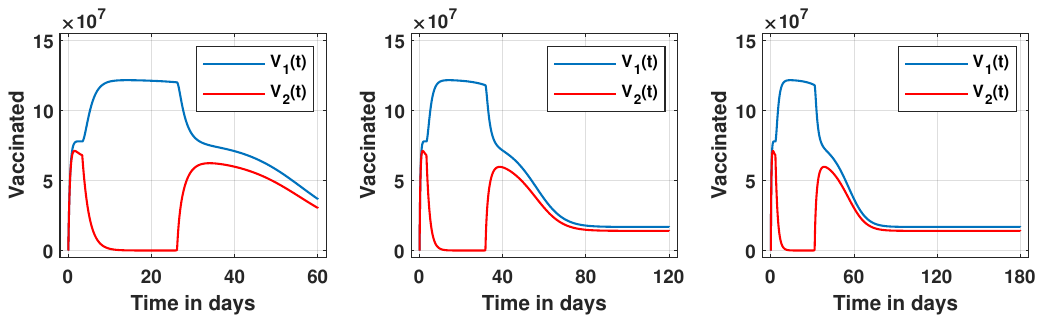}
\end{center}
\caption{\label{vaccinatedsensibility} Vaccinated for $\alpha_1 = 0.08, \alpha_2 = 0.064, \epsilon_1 = 0.432, \epsilon_2 = 0.54$.}
\end{figure}

For the control shown in Figure \ref{sensibility3}, the curve of vaccinated individuals is presented in Figure \ref{vaccinatedsensibility}. We also calculate the integral of the curves of vaccinated to observe the quantity of vaccines used during the campaign, as detailed in Table \ref{purchasedistribution7}.

\begin{table}[H]
\caption{\label{purchasedistribution7} Distribution for vaccine procurement according to Figure \ref{sensibility3}.}
\begin{center}
\begin{tabular}{c|c|c|c}
\hline
\textbf{Vaccine} & \textbf{T = 60 days} & \textbf{T = 120 days} & \textbf{T = 180 days} \\
\hline
$V_1$ & 71,60 \% & 71,21\% & 68,37 \% \\	     
\hline
$V_2$ & 28,40 \% & 28,79\% & 31,63 \% \\	     
\hline
\end{tabular}
\end{center}
\end{table}

\begin{figure}[H]
\centering
\subfigure[Optimal control $u_2(t)$ for $\theta_2$ fixed and $\theta_1$ varying.]{
\includegraphics[width=0.47\textwidth, trim=55mm 80mm 55mm 80mm]{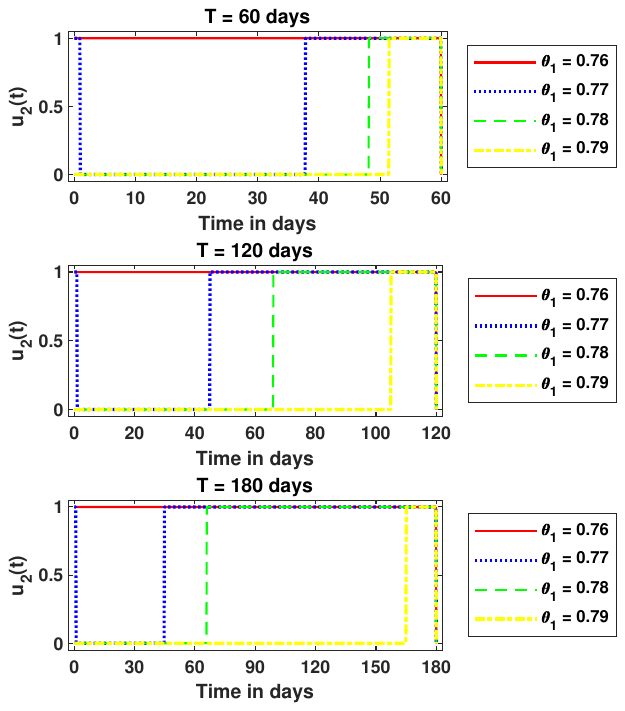}
\label{sensibility1}
}
\hfill
\subfigure[Optimal control $u_2(t)$ for $\theta_1$ fixed and $\theta_2$ varying.]{
\includegraphics[width=0.47\textwidth, trim=55mm 80mm 55mm 80mm]{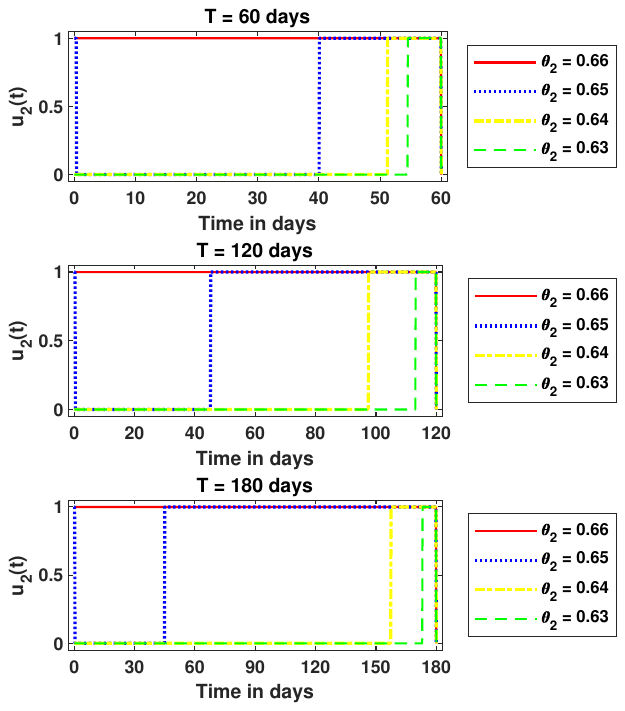}
\label{sensibility2}
}
\caption{Sensibility analysis for efficacy parameters $\theta_1$ and $\theta_2$.}
\label{fig:figuras}
\end{figure}

We also vary the efficacy parameters $\theta_1$ and $\theta_2$ to conduct a sensitivity analysis of these parameters. The behavior of the optimal control $u_1(t)$, as seen in Figure \ref{control6}, does not change, while the control $u_2(t)$ changes in the following way:
\begin{itemize}
\item Fixed $\theta_2=0.67$, the control $u_2(t)$ changes to values of $\theta_1$ above $\theta_1=0.77$, this is shown in Figure \ref{sensibility1};
\item Fixed $\theta_1=0.74$, the control $u_2(t)$ changes to values of $\theta_2$ below $\theta_2=0.65$, this is shown in Figure \ref{sensibility2};
\end{itemize}

In summary, as the efficacy parameters become closer, the simultaneous use of the vaccines is indicated. When the values diverge, the continued use of the $V_1$ vaccine, which has larger efficacy, is always indicated, while the use of the $V_2$ vaccine decreases.

\begin{figure}[H]
\begin{center}
\includegraphics[scale=0.8, trim=20mm 100mm 20mm 100mm]{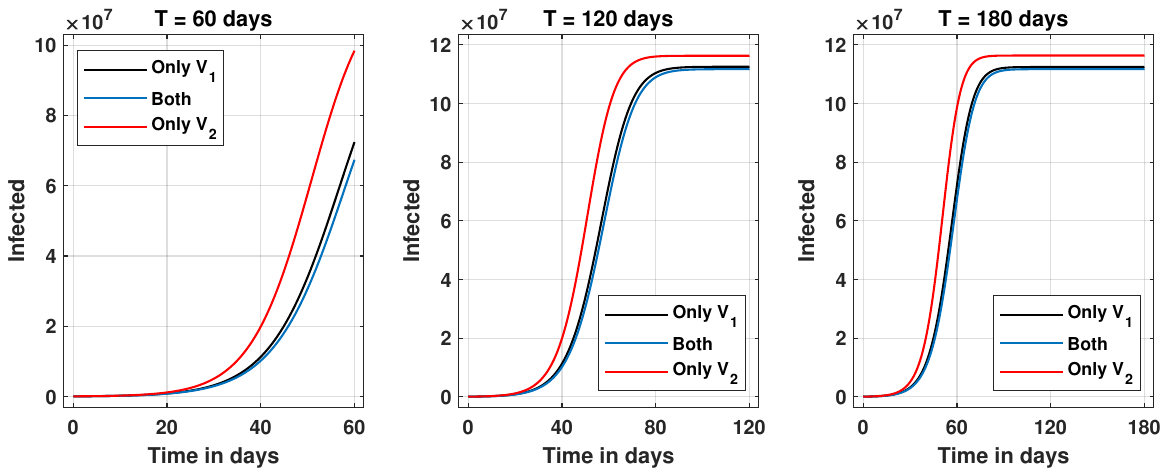}
\end{center}
\caption{\label{infected} Curves of infected for different use of the vaccines.}
\end{figure}

Finally, for vaccines $V_1$ with $74\%$ efficacy and $V_2$ with $67\%$ efficacy, we will see in Figure \ref{infected} the curves of infected people in three contexts: when we only use the $V_1$ vaccine, when we use both vaccines according to the optimal control in Figure \ref{control6}, and when we only use the $V_2$ vaccine.

Analyzing Figure \ref{infected}, we highlight the economic importance of optimal control. When we use only the $V_2$ vaccine, we have the lowest possible cost, but there is a considerable increase in the number of infected people compared to the other two strategies. Note that the infected curve exhibits practically the same behavior when we use both vaccines or when we use only the $V_1$ vaccine, which is more effective and has the most expensive implementation for the government. In this case, the strategy of using both vaccines, as defined by optimal control, presents a slight reduction in the number of infected people, with an affordable cost.

% ----------------------------------------------------------
\section{Conclusions}
% ----------------------------------------------------------

In this paper, we considered a SEIR model adapted for a vaccination campaign with two types of vaccines, $V_1$ and $V_2$. The controls of the system are the usage parameters of the vaccines. We found the optimal controls via Pontryagin’s maximum principle and used the forward-backward sweep method, to calculate the controls numerically.

Throughout the paper, we utilized four values for the efficacy of the vaccines: $91\%, 74\%, 67\%, 51\%$, and investigated the optimal controls for each pair of these values. Based on what was indicated by the optimal control, we showed how the procurement of vaccines should be distributed and how they should be used. When the efficacy values were very similar, we presented a sensitivity analysis for the parameters. Furthermore, we showed how the simultaneous use of the two vaccines, according to the optimal control, reduces the number of infected individuals at a more affordable cost.

When the $V_2$ vaccine is $51\%$ effective, the use of the other vaccine with efficacy greater than or equal to $67\%$ is indicated, i.e., a $51\%$ effective vaccine is only recommended when it is the only option available. We understand that this occurs due to the effect of efficacy on the transmission rate, which decreases linearly. The $V_2$ vaccine has a transmission rate that is considered high when compared to the rates of other vaccines. 

When considering the $V_1$ vaccine with $91\%$ efficacy and the $V_2$ vaccine with $74\%$ or $67\%$ efficacy, the use of the vaccine with higher efficacy is indicated almost all the time, followed by a simultaneous use of both vaccines during the last days of the campaign. This likely occurs to balance economic expenses, as the population can achieve collective immunity. The distribution for the purchase of vaccines indicated that in general more than $93\%$ of the vaccines utilized must be from the $V_1$ vaccine.

When the $V_1$ vaccine is $74\%$ effective and the $V_2$ vaccine is $67\%$ effective, the optimal control indicated a simultaneous use throughout the entire campaign. In this case, the procurement of both vaccines is recommended. In the sensitivity analysis of the recovery rate and the rate of returning to the susceptible class, the control of simultaneous use only changes when the recovery rate of the $V_2$ vaccine and the rate of returning to the susceptible class of the $V_1$ vaccine decrease by $20\%$. In this scenario, for all campaigns lengths, the exclusive use of the vaccine with higher efficacy is indicated for approximately the first 30 days, followed by simultaneous use until the end of the campaign. The purchase should consist of approximately $70\%$ of the $V_1$ vaccine and $30\%$ of the $V_2$ vaccine. For the sensitivity of efficacy values, when the $V_2$ vaccine is fixed with efficacy of $67\%$, the scenario of simultaneous use only changes when the $V_1$ vaccine has efficacy above $77\%$. Conversely, when the $V_1$ vaccine is fixed with $74\%$ of efficacy, the scenario of simultaneous use is only not indicated when the $V_2$ vaccine has efficacy lower than $65\%$. When these changes occurs, it is recommended to use only the $V_1$ vaccine for an initial period, followed by simultaneous use.

Finally, we also observed that the context in which simultaneous use is indicated by optimal control throughout the entire campaign presents a smaller number of infected individuals when compared with the contexts in which only one vaccine is used, and this minimizes economic costs.

We conclude that unless we have a vaccine with very low effectiveness (e.g., $51\%$, for which use is not indicated) or two vaccines with very similar efficacy (for which simultaneous use is indicated at all times), the optimal controls indicate a standard behavior: initial use of the more effective vaccine only, followed by simultaneous use.

In our model, we utilized COVID-19 data from Brazil for the initial conditions and data from the literature for parameter values. This type of study can be conducted to investigate vaccination campaigns for other diseases by adapting the model according to their characteristics and using initial conditions that suit the specific context (for example, Brazil is very large, and there may be underreporting and other factors that introduce noise in the data), followed by appropriate estimation of the parameters.

\section*{Acknowledgements}
The first author was partially supported by CAPES, Brazil.

%\bibliographystyle{unsrt}
%\bibliography{references}  

\end{document}